\documentstyle[12pt,amssymb]{article}
\title{On a function of Marcel Riesz}
\author{Gene Ward Smith}
\date{}

\newtheorem{theorem}{Theorem}

\newtheorem{corollary}[theorem]{Corollary}

\def\Riesz{\mathop{\rm Riesz}}

\begin{document}

\maketitle

\begin{abstract}

Formulas for calculating the Riesz function, introduced by Marcel Riesz
in connection with the Riemann hypothesis, are derived; and the behavior 
of the Riesz function is discussed.

\end{abstract}

In a paper of 1916 (\cite {R}), Marcel Riesz introduced an entire function which we 
will call the Riesz function. It may be defined by its Maclaurin series,
which is
$$\Riesz(x) = \sum_{k=1}^\infty \frac{(-1)^{k+1}x^k}{(k-1)! \zeta(2k)}$$
The interest in the function is that $\Riesz(x) = O(x^{\frac14 + \epsilon})$
is equivalent to the Riemann hypothesis.

Some interest, therefore, attaches to computing its values for positive real
values of x. We can convert the Maclaurin series into one with rational
coefficients by using the coefficients of the Laurent series expansion 
of the inverse hyperbolic cotangent. If
$$\frac{x}{2} \coth \frac{x}{2} = \sum_{n=0}^\infty c_n x^n = 
\sum_{n=0}^\infty \frac{B_{2n}x^{2n}}{(2n)!} =
1 + \frac{1}{12} x^2 - \frac{1}{720}x^4 + \cdots$$
where $B_n$ is the nth Bernoulli number, then using 
$\zeta(2k) = (-1)^{k+1}{\frac12} B_{2k} (2 \pi)^{2k}/(2k)!$, we find
$$\frac12 \Riesz(4 \pi^2 x) = \sum_{n=1}^\infty \frac{x^n}{c_{2n}(n-1)!} = 12x - 720x^2 + 15120x^3 - \cdots$$
As might be expected, this is not an efficient way of computing the Riesz function
unless x is small; however we can transform it into another infinite
series with much better convergence properties.

\begin{theorem}
Let 
$$\sum_{n=1}^\infty \frac{a_n}{n^s}$$
be a Dirichlet series absolutely convergent in the
region $\Re(x) > \sigma$, and let $f(s)$ be the function defined by its sum. Then if $c>\sigma$
is a positive constant,
let $R_c$ be the (entire) function defined by the power series
$$R_c(x) = \sum_{k=1}^\infty (-1)^k \frac{f(ck)x^k}{(k-1)!}$$
Then the series
$$x \sum_{n=1}^\infty \frac{a_n}{n^c} \exp(- \frac{x}{n^c})$$
converges absolutely for all values of x to $R_c(x)$.
\end{theorem}

Proof:
$f(s)$ is bounded on the region $s \ge c$, and hence $f(ck)$ is bounded for all positive integers 
k. Hence the power series converges everywhere by comparison with the series for $x \exp(-x)$,
defining an entire function.
Expanding $f(ck)$ in the power series expression gives
$$R_c(x) = \sum_{k=1}^\infty \frac{(-1)^{k+1}x^k}{(k-1)!}( \sum_{n=1}^\infty a_n n^{-ck})$$
Because of absolute convergence, we may expand this out to a double series and reverse 
the order of summation
$$R_c(x) = \sum_{k=1}^\infty \sum_{n=1}^\infty a_n \frac{(-1)^{k+1}(x/n^c)^k}{(k-1)!} = 
\sum_{n=1}^\infty a_n \sum_{k=1}^\infty \frac{(-1)^{k+1}(x/n^c)^k}{(k-1)!}$$
Collecting the terms in k now gives the theorem.
$\Box$

We may improve the convergence properties of the above series by using Kummer's method
of accelerated convergence.

\begin{corollary}
Under the conditions of the theorem, we have
$$R_c(x) = \sum_{k=1}^{m} (-1)^{k+1} \frac{f(ck)x^k}{(k-1)!}
 + x( \sum_{n=1}^\infty \frac{a_n}{n^c} F_m(-\frac{x}{n^c}))$$
where $F_m(x)$ is the mth order remainder term in the Maclaurin series
expansion of $\exp(x)$:
$$F_m(x) = \exp(x) - \sum_{k=0}^{m-1}\frac{x^k}{k!}$$
\end{corollary}
Proof:
Expand out
$$R_c(x) = \sum_{k=1}^{m} (-1)^{k+1} \frac{f(ck)x^k}{(k-1)!} + x( \sum_{n=1}^\infty \frac{a_n}{n^c} (\exp(-\frac{x}{n^c}) - 
\sum_{k=0}^{m-1} (-1)^k \frac{x^k}{n^{ck}k!}))$$
$\Box$

Because the above series for $R_c$ converges uniformly, it may be differentiated term by term.

\begin{corollary}
For all values of x, we have
$$R_c'(x) = \frac{R_c(x)}{x} - R_{2c}(x)$$ 
\end{corollary}

Applying these more general results to the case where $f(s) = (\zeta(s))^{-1}$ and $c=2$ gives
\begin{corollary}
For all values of x, the series
$$x \sum_{n=1}^\infty \frac{\mu(n)}{n^2} \exp(-\frac{x}{n^2})$$
converges absolutely to $\Riesz(x)$, where $\mu(n)$ is the M\"{o}bius function.
\end{corollary}

\begin{corollary}
For all values of x, we have
$$\Riesz(x) = 2\sum_{k=1}^m \frac{(-1)^{k+1}(x/4\pi^2)^k}{c_{2k}(k-1)!}+
x( \sum_{n=1}^\infty \frac{\mu(n)}{n^2} F_m(-\frac{x}{n^2}))$$
where $c_k$ are the coefficients of $\frac{x}{2} \coth \frac{x}{2}$.
\end{corollary}

By Cauchy's form of Taylor's theorem, we have
$F_m(x) = (-1)^m \exp(-\xi) x^m/m!$
for $0 < \xi < x$. From 
this we may conclude that $|F_m(x)| < x^m/m!$, and hence
$$x\ |\frac{F_m(-x/n^2)}{n^2}| < \frac{x^{m+1}}{n^{2+2m}m!}$$
Fixing x and m, the terms are $O(n^{-2-2m})$. If we set n to be 
${\rm ceil}(\sqrt{x})$, the bound becomes less than
$1/m!$, so choosing an m which makes this less than some
error limit allows us to find the Riesz function for positive
real numbers to a reasonably well specified tolerance.

In order to explore the connection between the Riesz function and
the Riemann zeta function, we need its Mellin transform.

\begin{theorem}
Let $\sigma$ be the abscissa of absolute convergence of the Dirichlet
series for $f(s)$, that is, the value beyond which the series
$\sum |a_n|n^{-s}$ converges. Then if 
$-1 < \Re(s) < -\sigma/c$, we have
$$\int_0^\infty R_c(x)x^s \frac{dx}{x} = f(-cs)\Pi(s)$$
\end{theorem}

Proof:
The function being integrated is a sum of functions 
$$\sum_{n=1}^\infty \frac{a_n}{n^c}x^s \exp(-\frac{x}{n^c})$$
Term by term integration of a series of positive terms is valid if
either side converges. Since
$$\int_0^\infty \frac{|a_n|}{n^c}x \exp(-\frac{x}{n^c}) x^s \frac{dx}{x} = |a_n| n^{cs} \Pi(s)$$ 
the Mellin transform of its terms may be taken so long as $s>-1$ and
$\sum_{n=1}^\infty |a_n| n^{cs}$ converges, which it will iff
$s < - \sigma/c$. By dominated convergence we then have that
term by term integration is valid for the original series, giving
the result.

\begin{theorem}
The Riemann hypothesis is equivalent to the claim
$${\rm Riesz}_c(x) = O(x^{\frac{1}{2c}+\epsilon})$$
for any value of $c>1$.
\end{theorem}

Proof:
The inverse Mellin transform tells us that
$${\rm Riesz}_c(x) = \frac{1}{2 \pi i} \int_{\gamma-i \infty}^{\gamma+i \infty} 
x^{-s} \frac{\Pi(s)}{\zeta(-cs)}ds$$
where $-1 < \gamma < -1/c$. The strip of analyticity of a 
Mellin transform is the widest region $a < \Re(s) < b$ in which
the function is analytic, and in terms of this strip, the original has
growth $O(x^{-a+\epsilon})$ as x tends to zero, and $O(x^{-b+\epsilon})$
as x tends to infinity. The Riemann hypothesis is equivalent to the 
claim that the strip of analyticity is $-1 < \Re(s) < -\frac{1}{2c}$. 
In terms of the original of the transform, this means ${\rm Riesz}_c(x)$ is
$O(x^{1+\epsilon})$ as x tends towards zero, which is trivial, and
${\rm Riesz}_c(x)$ must be $O(x^{\frac{1}{2c}+\epsilon})$
as x tends towards infinity, which gives us the theorem.

If we set
$$\kappa(s) = \frac{\Gamma(s/2)}{\pi^{s/2}}\zeta(s)$$
then the functional equation assumes the simple form 
$\kappa(s) = \kappa(1-s)$. 
In terms of this function, the inverse Mellin transform expression
for the Riesz function becomes
$$\Riesz(x) = \frac{1}{2 \pi i} \int_{\gamma-i \infty}^{\gamma+i \infty} 
x^{-s} \frac{\pi^{s+1}}{\sin (\pi s) \kappa(-2s)}ds$$
The function $\kappa(s)$ has poles at zero and one, and zeros only in the critical
strip; hence $\sin (\pi s) \kappa(-2s)$ has zeros at all integer values
other than zero, and in the critical strip $-1/2 < \Re(s) < 0$. If we move the line of integration from 
$-1 < \Re(s) < -1/2$ leftward, and evaluate resides, we obtain
$$\Riesz(x) = \sum_{k=1}^\infty \frac{(-1)^{k+1}(x/\pi)^k}{\kappa(2k)}$$
which is another form of the Maclaurin series defining the Riesz function.

If we move the line of integration to the right, the result is more
interesting. If $\rho$ is a simple nontrivial zero the Riemann zeta 
function, the residue of 
$$x^{-s} \frac{\pi^{s+1}}{\sin (\pi s) \kappa(-2s)}$$
at $-\rho/2$ is 
$$\frac{\pi\csc(\pi\rho/2)}{2 \kappa'(\rho)}(\frac{x}{\pi})^{\rho/2}
= \frac{\Pi(-\rho/2)}{2\zeta'(\rho)}x^{\rho/2}$$

At positive integers $n$, the residue is 
$$\frac{(-1)^{n}(\pi/x)^n}{\kappa(-2n)}$$
Since $\kappa(-2n) = \kappa(2n+1)$ this may also be written
$$\frac{(-1)^{n}(\pi/x)^n}{\kappa(2n+1)}$$
If we set
$$G(x) = \sum_{k=1}^\infty \frac{(-1)^{k+1}(x/\pi)^k}{\kappa(2k+1)}
= \pi^{-1/2}\sum_{k=1}^\infty \frac{(-1)^{k+1}x^k}{(k-1/2)!\zeta(2k+1)}$$
then by comparison with the series for the Riesz function we see
that $G$ is entire. Assuming for convenience that all of the zeros of the
Riemann zeta function are simple, and leaving aside the
question of convergence, we have
$$\Riesz(x) = \sum_\rho \frac{\Pi(-\rho/2)}{2 \zeta'(\rho)}x^{\rho/2}
- G(\frac{1}{x})$$

The question of the convergence of this series is an interesting one. Hardy and Littlewood (\cite {HL})
considered a very similar sum, and their method applies to this one. Except for a proportion 
of imaginary parts of s which can be made as small as we like, we can justify moving across
the singularities and evaluating the resides, so that the series converges conditionally
if we group together zeros which are proportionally (compared to the $\log t$ density of zeros)
close together in terms of imaginary part. With a properly chosen constant $K$, we can produce 
conditional convergence by bracketing terms where if $g_1$ and $g_2$ are two imaginary parts of
a zeros of the zeta function, then we put them in the same bracket if
$$|g_1 - g_2| < K^{-g_1/\log g_1} + K^{-g_2/\log g_2},$$
for some constant $K > 1$. Since on average the imaginary parts of contiguous zeros
are $1/\log g$ apart, this brackets zeros only when they are close together.

In fact, however, the series almost certainly
converges absolutely. If we assume the Riemann hypothesis, the numerators $(-\rho/2)!$ will
be $\Gamma(\frac34 + i \frac{t}{2})$, which is $O(t^{\frac14}\exp(-\frac{\pi}{4}t))$. If
$1/\zeta'(\rho)$ is of polynomial growth, the series will converge absolutely, and the
standard conjecture about this rate of growth, due to Gonek, is that it is $O(t^{\frac13})$; this
is supported both by Grand Unified Ensemble considerations and the numerical evidence. Adding an
extra factor of $\log t$ because there are $O(\log t)$, rather than one, zeros between
successive integer values of t, and the sum should have the convergence properties of
a series of terms which are $O(n^{\frac{7}{12}} \log n \exp(-\frac{\pi}{4} n))$, and therefore
converge absolutely.

\end{document}